\documentclass[preprint]{imsart}

\RequirePackage[authoryear]{natbib}
\usepackage{amsmath,amssymb,amsfonts,amsthm}
\usepackage[colorlinks,citecolor=blue,linkcolor=blue,urlcolor=blue]{hyperref}
\usepackage{booktabs}
\usepackage{tikz}
\usetikzlibrary{arrows.meta}
\usepackage{microtype}

\startlocaldefs
\numberwithin{equation}{section}
\newcommand{\N}{\mathbb{N}}
\newcommand{\PP}{\mathbb{P}}
\newcommand{\E}{\mathbb{E}}
\newcommand{\cA}{\mathcal{A}}
\newcommand{\cB}{\mathcal{B}}
\newcommand{\Def}{\mathcal{D}}
\newcommand{\bn}{\mathbf{n}}
\newcommand{\rf}[2]{(#1)_{#2}}
\DeclareMathOperator{\Beta}{Beta}
\newcommand{\repourl}{\url{https://github.com/Raydonal/gibbs-allelic-rigidity}}

\theoremstyle{plain}
\newtheorem{thm}{Theorem}[section]
\newtheorem{prop}[thm]{Proposition}
\newtheorem{lem}[thm]{Lemma}
\newtheorem{cor}[thm]{Corollary}
\theoremstyle{definition}
\newtheorem{defn}[thm]{Definition}
\endlocaldefs

\begin{document}

\begin{frontmatter}

\title{Allelic partitions of \protect\texorpdfstring{$\Lambda$}{Lambda}-coalescents given the
number of alleles: sufficiency and Gibbs rigidity}
\runtitle{Allelic partitions given $K_n$}
\runauthor{R. R. Albuquerque, A. M. S. Mac\^edo and R. Ospina}

\begin{aug}
\author[A,B]{\fnms{Roberta R.}~\snm{Albuquerque}},
\author[C]{\fnms{Ant\^onio M. S.}~\snm{Mac\^edo}}
\and
\author[A,D]{\fnms{Raydonal}~\snm{Ospina}}
\address[A]{Departamento de Estat\'istica, Universidade Federal de Pernambuco,
Recife 50670-901, PE, Brazil.}
\address[B]{Escola Polit\'ecnica de Pernambuco, Universidade de Pernambuco,
Recife 50720-001, PE, Brazil.\\ \texttt{roberta.ralbuquerque@upe.br}}
\address[C]{Laborat\'orio de F\'isica Te\'orica e Computacional, Departamento de F\'isica,
Universidade Federal de Pernambuco, Recife 50670-901, PE, Brazil.\\ \texttt{antonio.smacedo@ufpe.br}}
\address[D]{Departamento de Estat\'istica, Universidade Federal da Bahia,
Salvador 40170-110, BA, Brazil.\\ \texttt{raydonalmartinez@ufba.br} (corresponding author)}
\end{aug}

\begin{abstract}
Under Kingman's coalescent the allelic partition of a sample of size $n$, given its number
$K_n$ of alleles, does not depend on the mutation rate $\theta$ and has the weights of the
Ewens sampling formula. We ask which $\Lambda$-coalescents keep either property. Closed
forms for the configurations with at least $n-2$ alleles show that, when the first $n-1$
individuals form one pair and singletons, individual $n$ joins the pair with at least twice
the probability that it joins a given singleton, with equality only for Kingman's
coalescent. Consequently, for each $n\ge4$, $K_n$ is sufficient for $\theta$ if and only if
$\Lambda$ is concentrated at $0$ or at $1$; and for each fixed $\theta$ the partition of a
sample of size six is of Gibbs type if and only if $\Lambda$ is concentrated at $0$, while
size five does not suffice. For Gibbs types in $(0,1)$ the proof rests on a polynomial with
$530$ positive integer coefficients. We also give the limits of the associated edge index
for Beta and Bolthausen--Sznitman coalescents.
\end{abstract}

\begin{keyword}[class=MSC2020]
\kwd[Primary ]{60G09}
\kwd{60J90}
\kwd[; secondary ]{92D15}
\kwd{05A18}
\kwd{62B05}
\end{keyword}

\begin{keyword}
\kwd{\protect\texorpdfstring{$\Lambda$}{Lambda}-coalescent}
\kwd{infinitely-many-alleles model}
\kwd{exchangeable random partition}
\kwd{Gibbs-type partition}
\kwd{Ewens sampling formula}
\kwd{sufficient statistic}
\end{keyword}

\end{frontmatter}

\section{Introduction}\label{sec:intro}

Let a sample of $n$ individuals have the genealogy of a $\Lambda$-coalescent
\citep{Pitman99,Sagitov99}, where $\Lambda$ is a finite nonzero measure on $[0,1]$, and let
mutations occur along every lineage at rate $\theta/2$, each creating an allele not seen
before (the infinitely-many-alleles model). The \emph{allelic partition} $\Pi_n$ of
$[n]=\{1,\dots,n\}$ groups the individuals that carry the same allele, and $K_n$ denotes
its number of blocks. For Kingman's coalescent \citep{Kingman1982}, $\Lambda=c\delta_0$,
the partition $\Pi_n$ follows the Ewens sampling formula with parameter $\theta/c$
\citep{Ewens1972}. Two properties of this formula concern the conditional law of $\Pi_n$
given $K_n$: it does not depend on $\theta$, so that $K_n$ is a sufficient statistic for
$\theta$; and it gives each partition of $[n]$ with $k$ blocks of sizes $n_1,\dots,n_k$ a
probability proportional to $\prod_j(n_j-1)!$. This paper asks which $\Lambda$-coalescents
keep either property.

The law of $\Pi_n$ is determined by a recursion due to \citet{Mohle06}, which
\citet[eq.~(2.7)]{Dong07} write for the exchangeable partition probability function (EPPF);
we recall it in Proposition~\ref{prop:recursion}.
The partitions of Gibbs type \citep{Pitman05} generalise the second property: $\Pi_n$
is of Gibbs type $\sigma\in[-\infty,1)$ when, given $K_n=k$, each partition of $[n]$ with
$k$ blocks of sizes $n_1,\dots,n_k$ has a probability proportional to
$\prod_j\rf{1-\sigma}{n_j-1}$, where $\rf{x}{m}=x(x+1)\cdots(x+m-1)$; for $\sigma=-\infty$
the weight is $1$. The Ewens formula is of type $0$. \citet[Theorem~3.1]{Mohle07} showed that the
sampling formula of a $\Lambda$-coalescent is regenerative if and only if $\Lambda$ is
concentrated at $0$ or at $1$, and \citet[Proposition~3.2]{Mohle07} that it coincides,
for all sample sizes, with a member of Pitman's two-parameter family only in the Ewens
case. Corollary~\ref{cor:rigidity} extends the second statement from the two-parameter
family to all partitions of Gibbs type, and Theorem~\ref{thm:rigidity} shows that, for each
fixed $\theta$, the sample size six already decides it.

Our first result concerns the configurations of $\Pi_n$ with at least $n-2$ blocks, which
we call edge configurations because $K_n$ is then within two of its largest value $n$.
Theorem~\ref{thm:edge} gives the probabilities of the configurations $1^n$, $(2,1^{n-2})$,
$(3,1^{n-3})$ and $(2,2,1^{n-4})$ in closed form, for every $\Lambda$ and every $n$, in
terms of the moments $h_j=\int(1-x)^j\Lambda(dx)$. The ratio
\[
R_n=\frac{p_n(3,1^{n-3})}{p_n(2,2,1^{n-4})}
\]
of the probabilities of two partitions of $[n]$ with these configurations has a predictive
meaning: given that the first $n-1$ individuals form one pair and $n-3$ singletons, $R_n$
is the probability that individual $n$ joins the pair divided by the probability that it
joins a given singleton. Under the Ewens formula $R_n=2$, the ratio of the two block
sizes. We show that $R_n\ge2$ for every $\Lambda$ and every $n\ge4$, with equality only if
$\Lambda=c\delta_0$, and that $R_n$ is an affine and strictly increasing function of
$\theta$ unless $\Lambda$ is concentrated at $0$ or at $1$. A partition of Gibbs type
$\sigma\in(-\infty,1)$ has $R_n=(2-\sigma)/(1-\sigma)$; accordingly we call
$\sigma^{(n)}=(R_n-2)/(R_n-1)\in[0,1]$ the \emph{edge index}. Three consequences follow.

\emph{Sufficiency} (Theorem~\ref{thm:sufficiency}). For every $n\ge4$, $K_n$ is sufficient
for $\theta$ if and only if $\Lambda$ is concentrated at $0$ or at $1$. Given $K_n=n-2$,
the probability that the configuration is $(3,1^{n-3})$ rather than $(2,2,1^{n-4})$ is at
least $8/(3n-1)$, its value under the Ewens formula, with equality only for Kingman's
coalescent. If $\Pi_n$ is of the same Gibbs type for two values of $\theta$, then
$\Lambda=c\delta_0$.

\emph{Rigidity at a fixed $\theta$} (Theorem~\ref{thm:rigidity}). For each fixed
$\theta>0$, $\Pi_6$ is of Gibbs type, for some $\sigma\in[-\infty,1)$, if and only if
$\Lambda=c\delta_0$. For types $\sigma\le0$ this holds at $n=4$ already, by $R_4\ge2$. For
$\sigma\in(0,1)$, the identities $R_n=(2-\sigma)/(1-\sigma)$ at $n=4,5,6$ determine $h_2$,
$h_3$ and $h_4$ from $\sigma$, $\theta$, $h_0$ and $h_1$, and the Gibbs identity between the
configurations $(3,2,1)$ and $(2,2,2)$ then fails. This failure is the positivity of a
rational function whose numerator contains a polynomial in three variables of degree $18$
with $530$ integer coefficients, all positive (Lemma~\ref{lem:certificate}); the
coefficients are supplied as data. The sample size six cannot be lowered to five: for every
$\theta>0$, $\Pi_5$ is of Gibbs type $2/3$ when $\Lambda$ is $3\theta/2$ times the Lebesgue
measure on $[0,1]$, a time change of the Bolthausen--Sznitman coalescent
(Proposition~\ref{prop:sharp}).

\emph{Large samples} (Proposition~\ref{prop:limits}). As $n\to\infty$ with $\theta$ fixed,
$\sigma^{(n)}\to0$ if $\Lambda$ has an atom at $0$;
$\sigma^{(n)}\to(2-\alpha)(1+\alpha)/(1+3\alpha-\alpha^2)$, whatever $\theta$, if $\Lambda$
has a density asymptotic to $Cx^{1-\alpha}$ as $x\downarrow0$ with $1<\alpha<2$, as
$\Beta(2-\alpha,\alpha)$ has; $\sigma^{(n)}\to(4+\theta)/(6+\theta)$ for the
Bolthausen--Sznitman coalescent; and $\sigma^{(n)}\to1$ for such densities with
$0<\alpha<1$. The edge configurations have probabilities tending to $0$
(Lemma~\ref{lem:rare}), so these limits concern the law of $\Pi_n$ on rare events; for
$\Beta(2-\alpha,\alpha)$ with $1<\alpha<2$ we compare them with the limits of the block
counts obtained by \citet{BBS2007}, which describe typical samples.

\citet{Joyce1998} asked whether the Ewens family is the only one-parameter family of
partition structures in which the number of types is sufficient. By his Theorem~1, $K_n$ is
sufficient for every $n$ exactly when the predictive probability of a new type depends only
on the number of types and the parameter, and the choice among the existing types does not
depend on the parameter; Pitman's two-parameter family with one parameter held fixed is of
this kind, and so is his Example~1, in which $K_n$ determines the whole configuration. By
his Theorem~2, $K_n$ is not sufficient for infinite-alleles models with selection whose
representing measure has a density with respect to the Poisson--Dirichlet law, unless that
density is $1$. Theorem~\ref{thm:sufficiency} concerns neutral models with multiple-merger
genealogies, at each fixed sample size $n\ge4$: $K_n$ is sufficient only for Kingman's
coalescent and for the star-shaped coalescent, $\Lambda=c\delta_1$. The latter is
degenerate in the sense of Joyce's Example~1, since a partition with $k$ blocks then always
has the configuration $(n-k+1,1^{k-1})$.

The partitions studied here are allelic. Genealogical partitions behave differently: for
the Bolthausen--Sznitman coalescent the partition of $\N$ into the blocks present at a
fixed time $t$ has the Poisson--Dirichlet law $\mathrm{PD}(e^{-t},0)$ of Pitman's
two-parameter family \citep{Bolthausen98,Pitman99}; see \citet[Theorem~6.2]{Berestycki09}.

Section~\ref{sec:model} recalls the model and the recursion, and Section~\ref{sec:gibbs}
the partitions of Gibbs type on $[n]$. Section~\ref{sec:edge} treats the configurations
with at least $n-2$ blocks, Section~\ref{sec:sufficiency} sufficiency,
Section~\ref{sec:rigidity} the rigidity at a fixed $\theta$ and Section~\ref{sec:limits}
large samples. Section~\ref{sec:open} lists open problems, and
Appendix~\ref{app:certificate} describes the polynomials of Lemma~\ref{lem:certificate}.

\section{Allelic partitions of \texorpdfstring{$\Lambda$}{Lambda}-coalescents}\label{sec:model}

Throughout, $\Lambda$ is a finite nonzero measure on $[0,1]$ and $\theta>0$. Statements
are written in terms of $\theta$; in proofs we abbreviate $r=\theta/2$, the mutation rate
per lineage. We write $\binom{b}{m}$ for binomial coefficients,
$\rf{x}{m}=x(x+1)\cdots(x+m-1)$ with $\rf{x}{0}=1$, and $0^0=1$.

\begin{figure}[t]
\centering
\begin{tikzpicture}[
  font=\footnotesize,
  lin/.style={draw=black!70,line width=0.8pt},
  mut/.style={fill=red!75!black,draw=none},
  lab/.style={font=\footnotesize},
  blk/.style={rounded corners=3pt,draw=black!45,inner sep=3pt,font=\footnotesize}]

\foreach \i in {1,...,8} { \coordinate (L\i) at (\i*1.05,0); }
\node[lab] at (1*1.05,-0.32) {$1$};  \node[lab] at (2*1.05,-0.32) {$2$};
\node[lab] at (3*1.05,-0.32) {$3$};  \node[lab] at (4*1.05,-0.32) {$4$};
\node[lab] at (5*1.05,-0.32) {$5$};  \node[lab] at (6*1.05,-0.32) {$6$};
\node[lab] at (7*1.05,-0.32) {$7$};  \node[lab] at (8*1.05,-0.32) {$8$};

\def\ya{1.15}
\coordinate (A) at (3*1.05,\ya);
\draw[lin] (L2)--(2*1.05,\ya)--(4*1.05,\ya);
\draw[lin] (L3)--(3*1.05,\ya);
\draw[lin] (L4)--(4*1.05,\ya);
\def\yb{1.75}
\coordinate (B) at (6.5*1.05,\yb);
\draw[lin] (L6)--(6*1.05,\yb)--(7*1.05,\yb);
\draw[lin] (L7)--(7*1.05,\yb);
\def\yc{2.5}
\coordinate (C) at (2*1.05,\yc);
\draw[lin] (L1)--(1*1.05,\yc)--(3*1.05,\yc);
\draw[lin] (A)--(3*1.05,\yc);
\def\yd{3.1}
\coordinate (D) at (5.75*1.05,\yd);
\draw[lin] (L5)--(5*1.05,\yd)--(6.5*1.05,\yd);
\draw[lin] (B)--(6.5*1.05,\yd);
\def\ye{3.9}
\draw[lin] (C)--(2*1.05,\ye)--(8*1.05,\ye);
\draw[lin] (D)--(5.75*1.05,\ye);
\draw[lin] (L8)--(8*1.05,\ye);
\draw[lin] (4.9*1.05,\ye)--(4.9*1.05,\ye+0.45);

\def\mk#1{\draw[mut] (#1) circle (2.1pt);}
\mk{1*1.05,1.55}
\mk{3*1.05,2.05}
\mk{5*1.05,1.30}
\mk{6.5*1.05,2.15}
\mk{8*1.05,3.35}

\draw[-{Latex[length=2mm]},black!55] (0.28,0) -- (0.28,4.35);
\node[lab,rotate=90,anchor=south,black!55] at (0.06,2.1) {time into the past};
\node[lab,anchor=east,black!55] at (0.20,0) {present};

\node[blk,fill=blue!8]   at (1.0*1.05,-1.05) {$\{1\}$};
\node[blk,fill=green!10] at (3.0*1.05,-1.05) {$\{2,3,4\}$};
\node[blk,fill=orange!12]at (5.0*1.05,-1.05) {$\{5\}$};
\node[blk,fill=purple!8] at (6.6*1.05,-1.05) {$\{6,7\}$};
\node[blk,fill=black!6]  at (8.0*1.05,-1.05) {$\{8\}$};
\node[lab,anchor=east] at (0.55,-1.05) {$\Pi_8=$};
\end{tikzpicture}
\caption{A coalescent tree with mutations and the allelic partition it induces, for
$n=8$. Lineages run from the present upwards into the past and coalesce; a
$\Lambda$-coalescent allows more than two of them to merge at once, as in the triple merger
of $2,3,4$. Mutations (dots) fall on the branches at rate $\theta/2$, and each creates an
allele not seen before. Two individuals carry the same allele exactly when no mutation lies
between them and their common ancestor, so the first mutation met above a leaf closes its
class. These mutations occur at different times: the partition $\Pi_8$ at the bottom is not
a cross-section of the tree at one time.}
\label{fig:process}
\end{figure}
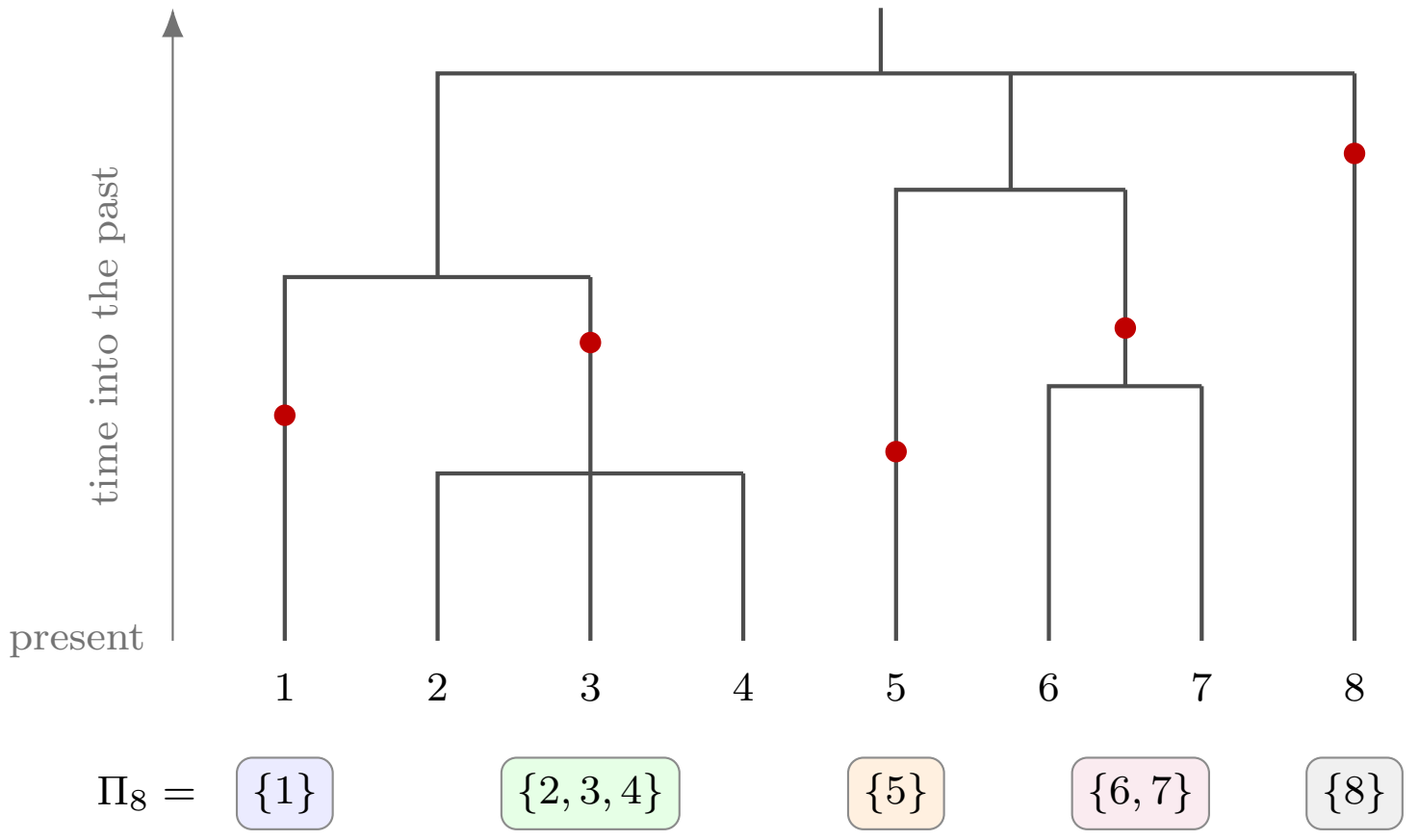

\paragraph{The model}
In the $\Lambda$-coalescent \citep{Pitman99,Sagitov99}, while there are $b$ blocks, each
particular set of $m$ of them, $2\le m\le b$, merges into one block at rate
\begin{equation}\label{eq:rates}
\lambda_{b,m}=\int_{[0,1]}x^{m-2}(1-x)^{b-m}\,\Lambda(dx).
\end{equation}
We follow the lineages of a sample of $n$ individuals backwards in time and superimpose on
each lineage the points of independent Poisson processes of rate $\theta/2$, the
mutations. Under the infinitely-many-alleles model every mutation creates a new allele, so
two individuals carry the same allele if and only if no mutation lies on the path that
joins them in the genealogical tree. Equivalently, a lineage is removed at its first
mutation, and the individuals of the sample that it carries at that moment form one allelic
class. The \emph{allelic partition} $\Pi_n$ is the partition of $[n]$ into these classes,
and $K_n$ is its number of blocks. The rates satisfy
$\lambda_{b,m}=\lambda_{b+1,m}+\lambda_{b+1,m+1}$, because $(1-x)+x=1$ under the integral,
so the process followed on the lineages of $[n]$ is the restriction of the process followed
on the lineages of $[n+1]$. Consequently $\Pi_n$ is exchangeable and is the restriction of
$\Pi_{n+1}$ to $[n]$, and $(\Pi_n)_{n\ge1}$ are the restrictions of an exchangeable random
partition $\Pi$ of $\N$. Figure~\ref{fig:process} shows an example.

For a composition $\bn=(n_1,\dots,n_k)$ of $n$, that is, positive integers with sum $n$,
let $p_n(\bn)$ be the probability that $\Pi_n$ equals a given partition of $[n]$ whose
blocks have the sizes $n_1,\dots,n_k$. By exchangeability, $p_n(\bn)$ depends neither on
the partition chosen nor on the order of the parts; $p_n$ is the EPPF of $\Pi_n$
\citep{Pitman02}. The \emph{configuration} of a partition is the integer partition formed
by its block sizes; $1^j$ stands for $j$ parts equal to $1$, $e_j$ for the $j$-th unit
vector, and $\bn^{(j)}$ for $\bn$ with its $j$-th part removed.

\paragraph{Moments}
Let $h_j=\int_{[0,1]}(1-x)^j\,\Lambda(dx)$ for $j\ge0$. The sequence $(h_j)$ is
nonnegative and nonincreasing, $h_0=\Lambda([0,1])>0$, and $h_j\to\Lambda(\{0\})$ as
$j\to\infty$ by dominated convergence.

\begin{lem}\label{lem:dictionary}
For $2\le m\le b$,
\begin{equation}\label{eq:dictionary}
\lambda_{b,m}=\sum_{i=0}^{m-2}(-1)^i\binom{m-2}{i}h_{b-m+i},
\end{equation}
in particular $\lambda_{b,2}=h_{b-2}$ and $\lambda_{b,3}=h_{b-3}-h_{b-2}$, and the total
merger rate with $b$ blocks is
\begin{equation}\label{eq:g}
g_b:=\sum_{m=2}^{b}\binom{b}{m}\lambda_{b,m}=\sum_{j=0}^{b-2}(j+1)h_j .
\end{equation}
Moreover $h_0-h_1=\int x\,\Lambda(dx)$ vanishes if and only if $\Lambda=h_0\delta_0$, and
$h_2=0$ if and only if $h_1=0$ if and only if $\Lambda=h_0\delta_1$.
\end{lem}

\begin{proof}
Formula \eqref{eq:dictionary} follows from $x^{m-2}=(1-(1-x))^{m-2}$. For \eqref{eq:g}, put
$y=1-x$. For $0\le y<1$,
\begin{align*}
\sum_{m=2}^{b}\binom{b}{m}x^{m-2}y^{b-m}&=\frac{1-y^b-b(1-y)y^{b-1}}{(1-y)^2}\\
&=\frac{1-by^{b-1}+(b-1)y^b}{(1-y)^2}=\sum_{j=0}^{b-2}(j+1)y^j,
\end{align*}
and both ends equal $\binom{b}{2}$ at $y=1$; integrating against $\Lambda$ gives
\eqref{eq:g}. The last assertions hold because $x>0$ on $(0,1]$ and $(1-x)^j>0$ on
$[0,1)$.
\end{proof}

We write $\rho_b=g_b+b\theta/2$ for the total rate of events while there are $b$ lineages.
Replacing $(\Lambda,\theta)$ by $(c\Lambda,c\theta)$, $c>0$, multiplies every rate by $c$,
replaces each $h_j$ by $ch_j$, and leaves the law of $\Pi_n$ unchanged. For Kingman's
coalescent, $\Lambda=c\delta_0$, $h_j=c$ for all $j$, $\lambda_{b,2}=c$ and
$\lambda_{b,m}=0$ for $m\ge3$. For the star-shaped coalescent, $\Lambda=c\delta_1$,
$h_0=c$, $h_j=0$ for $j\ge1$ and $\lambda_{b,m}=c\,\mathbf 1\{m=b\}$.

\begin{prop}\label{prop:recursion}
$p_1(1)=1$, and for $n\ge2$ and every composition $\bn=(n_1,\dots,n_k)$ of $n$,
\begin{equation}\label{eq:recursion}
\rho_n\,p_n(\bn)=\frac\theta2\sum_{j:\,n_j=1}p_{n-1}(\bn^{(j)})
+\sum_{m=2}^{n}\lambda_{n,m}\sum_{j:\,n_j\ge m}\binom{n_j}{m}\,
p_{n-m+1}\bigl(\bn-(m-1)e_j\bigr).
\end{equation}
In particular $p_n$ is determined by $\theta$ and $h_0,\dots,h_{n-2}$.
\end{prop}

This is the recursion of \citet{Mohle06} in the form given by \citet[eq.~(2.7)]{Dong07} for
the EPPF; there, both sides are divided by $\rho_n$.

\begin{proof}
Fix a partition $\pi$ of $[n]$ with block sizes $\bn$, and consider the first event,
backwards in time, among the $n$ lineages of the sample. It is a mutation on a given
lineage with probability $r/\rho_n$, and the merger of a given set $S$ of $m\ge2$ lineages
with probability $\lambda_{n,m}/\rho_n$. If the lineage of individual $\ell$ mutates first,
no other individual carries its allele, so $\Pi_n=\pi$ requires $\{\ell\}$ to be a block of
$\pi$, say the $j$-th; the remaining $n-1$ lineages then evolve as the process started from
$n-1$ lineages, and the conditional probability that $\Pi_n=\pi$ is
$p_{n-1}(\bn^{(j)})$. If the lineages in $S$ merge first, the individuals in $S$ carry
the allele of the merged lineage, so $\Pi_n=\pi$ requires $S$ to lie in a block of $\pi$,
say the $j$-th, with $n_j\ge m$; the $n-m+1$ remaining lineages then evolve as the
process started from $n-m+1$ lineages, and the conditional probability is
$p_{n-m+1}(\bn-(m-1)e_j)$. The $j$-th block contains $\binom{n_j}{m}$ sets of
size $m$. Summing over the first event gives \eqref{eq:recursion}. Since $\rho_n>0$,
\eqref{eq:recursion} determines $p_n$ from $p_1,\dots,p_{n-1}$, and by
\eqref{eq:dictionary} and \eqref{eq:g} only $h_0,\dots,h_{n-2}$ enter.
\end{proof}

\begin{cor}\label{cor:ewens}
If $\Lambda=c\delta_0$, then for every $n\ge1$ and every composition $\bn$ of $n$ with
$k$ parts, $p_n(\bn)=\vartheta^k\prod_{j=1}^k(n_j-1)!/\rf{\vartheta}{n}$ with
$\vartheta=\theta/c$; that is, $\Pi_n$ follows the Ewens sampling formula with parameter
$\theta/c$.
\end{cor}

\begin{proof}
Let $E_n(\bn)$ denote the right-hand side; $E_1(1)=1$. Here
$\rho_n=c\binom n2+nr=\tfrac{cn}{2}(n-1+\vartheta)$, and only the term $m=2$ of
\eqref{eq:recursion} survives. For $n_j=1$,
$E_{n-1}(\bn^{(j)})=E_n(\bn)(\vartheta+n-1)/\vartheta$, and for $n_j\ge2$,
$E_{n-1}(\bn-e_j)=E_n(\bn)(\vartheta+n-1)/(n_j-1)$. Since $r/\vartheta=c/2$
and $\binom{n_j}{2}/(n_j-1)=n_j/2$, the right-hand side of
\eqref{eq:recursion} with $p=E$ equals
$E_n(\bn)(\vartheta+n-1)\tfrac c2\sum_jn_j=\rho_nE_n(\bn)$. By uniqueness of
the solution of \eqref{eq:recursion}, $p_n=E_n$.
\end{proof}

\section{Partitions of Gibbs type on \texorpdfstring{$[n]$}{[n]}}\label{sec:gibbs}

For $\sigma\in[-\infty,1)$ and a composition $\bn=(n_1,\dots,n_k)$ put
$w_\sigma(\bn)=\prod_{j=1}^k\rf{1-\sigma}{n_j-1}$ if $\sigma>-\infty$, and
$w_{-\infty}(\bn)=1$. For any exchangeable random partition of $[n]$, with EPPF $p_n$, and
any $\sigma$ we write $u_n(\bn)=p_n(\bn)/w_\sigma(\bn)$; the value of $\sigma$ will always
be clear from the context.

\begin{defn}\label{def:gibbs}
Let $\sigma\in[-\infty,1)$. An exchangeable random partition of $[n]$ is \emph{of Gibbs
type $\sigma$} if there are numbers $V_{n,k}$, $1\le k\le n$, such that
$p_n(\bn)=V_{n,k}\,w_\sigma(\bn)$ for every composition $\bn$ of $n$ with $k$ parts, that
is, if $u_n$ takes the same value $V_{n,k}$ on all compositions of $n$ with $k$ parts. An
exchangeable random partition of $\N$ is of Gibbs type $\sigma$ if its restriction to $[n]$
is, for every $n$.
\end{defn}

Since $w_\sigma>0$, the $V_{n,k}$ are nonnegative, and the definition says that for every
$k$ with $\PP(K_n=k)>0$ the conditional law of $\Pi_n$ given $K_n=k$ is proportional to
$w_\sigma$ of the block sizes. The weights and the recursion of Lemma~\ref{lem:addition}
below are those of the Gibbs partitions of \citet[Lemma~2 and Definition~3]{Pitman05}.
For $n\le3$ each value of $k$ corresponds to a single configuration, so every exchangeable
partition of $[n]$ is of Gibbs type $\sigma$ for every $\sigma$. For $n=4$ only $k=2$
carries two configurations, $(3,1)$ and $(2,2)$, so a partition of $[4]$ is of Gibbs type
$\sigma>-\infty$ if and only if $(1-\sigma)p_4(3,1)=(2-\sigma)p_4(2,2)$.

\begin{lem}\label{lem:addition}
Let $(\Pi_n)_{n\ge1}$ be the restrictions of an exchangeable random partition of $\N$, let
$\sigma\in[-\infty,1)$ and $n\ge2$. If $\Pi_n$ is of Gibbs type $\sigma$, then so is $\Pi_m$
for every $m\le n$, and for $1\le k\le m<n$
\begin{align*}
V_{m,k}&=(m-\sigma k)V_{m+1,k}+V_{m+1,k+1}&&(\sigma>-\infty),\\
V_{m,k}&=kV_{m+1,k}+V_{m+1,k+1}&&(\sigma=-\infty).
\end{align*}
\end{lem}

\begin{proof}
It suffices to take $m=n-1$. By consistency, for a composition $\bn$ of $m$ with $k$
parts, $p_m(\bn)=\sum_{j=1}^kp_{m+1}(\bn+e_j)+p_{m+1}(\bn,1)$: the element
$m+1$ joins the $j$-th block or forms a block of its own. Since
$\rf{1-\sigma}{n_j}=\rf{1-\sigma}{n_j-1}(n_j-\sigma)$, the Gibbs form at
$m+1$ gives
$p_m(\bn)=w_\sigma(\bn)\bigl[V_{m+1,k}\sum_j(n_j-\sigma)+V_{m+1,k+1}\bigr]
=w_\sigma(\bn)\bigl[(m-\sigma k)V_{m+1,k}+V_{m+1,k+1}\bigr]$, and for $\sigma=-\infty$
the sum $\sum_j1$ equals $k$.
\end{proof}

\begin{lem}\label{lem:predictive}
Let $\sigma\in[-\infty,1)$ and $n\ge2$. An exchangeable random partition of $[n]$ is of
Gibbs type $\sigma$ if and only if, for every composition $\bn$ of $n-1$ with $k$
parts and all $i\ne j$ in $\{1,\dots,k\}$,
\begin{equation}\label{eq:predictive}
\begin{aligned}
(n_j-\sigma)\,p_n(\bn+e_i)&=(n_i-\sigma)\,p_n(\bn+e_j)&&(\sigma>-\infty),\\
p_n(\bn+e_i)&=p_n(\bn+e_j)&&(\sigma=-\infty).
\end{aligned}
\end{equation}
\end{lem}

This is the familiar predictive form of Gibbs partitions; we include the proof for the
finite-$n$ version used below.
In predictive terms: if $\pi$ is a partition of $[n-1]$ and $B,B'$ are two of its blocks,
then $(|B'|-\sigma)\,\PP(n\in B\mid\Pi_{n-1}=\pi)=(|B|-\sigma)\,\PP(n\in B'\mid\Pi_{n-1}=\pi)$
whenever $\PP(\Pi_{n-1}=\pi)>0$; the probability that $n$ starts a new block is not
constrained.

\begin{proof}
If $p_n=V_{n,k}w_\sigma$, both sides of \eqref{eq:predictive} equal
$V_{n,k}w_\sigma(\bn)(n_i-\sigma)(n_j-\sigma)$. Conversely,
\eqref{eq:predictive} says $u_n(\bn+e_i)=u_n(\bn+e_j)$: $u_n$ takes the same value
on two compositions with $k$ parts when one is obtained from the other by moving a unit
from a part of size at least $2$ to another part. Moving units into the first part
connects every composition of $n$ with $k$ parts to $(n-k+1,1^{k-1})$, so $u_n$ is constant
on the compositions with $k$ parts.
\end{proof}

\section{Configurations with at least \texorpdfstring{$n-2$}{n-2} blocks}\label{sec:edge}

We call $1^n$, $(2,1^{n-2})$, $(3,1^{n-3})$ and $(2,2,1^{n-4})$ the edge configurations:
they are the configurations with at least $n-2$ blocks, for which $K_n$ is within two of
its largest value $n$.

\begin{thm}\label{thm:edge}
Let $n\ge4$ and
\begin{equation}\label{eq:AB}
\begin{aligned}
\cA_n&=\sum_{m=3}^{n}\bigl[3h_{m-2}g_{m-1}+(m-2)(h_{m-3}-h_{m-2})\rho_{m-1}\bigr],\\
\cB_n&=2\sum_{m=4}^{n}(m-3)\,h_{m-2}\,g_{m-1}.
\end{aligned}
\end{equation}
\begin{itemize}
\item[(a)] The edge configurations have the probabilities
\begin{align*}
p_n(1^n)&=\prod_{m=2}^{n}\frac{m\theta}{2\rho_m},\\
p_n(2,1^{n-2})&=\frac{2\,p_n(1^n)\,g_n}{n(n-1)\theta},\\
p_n(3,1^{n-3})&=\frac{4\,p_n(1^n)\,\cA_n}{n(n-1)(n-2)\theta^2},\\
p_n(2,2,1^{n-4})&=\frac{4\,p_n(1^n)\,\cB_n}{n(n-1)(n-2)(n-3)\theta^2}.
\end{align*}
Moreover $\cA_n>0$, so $R_n=p_n(3,1^{n-3})/p_n(2,2,1^{n-4})=(n-3)\cA_n/\cB_n$, with
$R_n=+\infty$ when $\cB_n=0$.
\item[(b)] $R_n\ge2$, with equality if and only if $\Lambda=c\delta_0$ for some $c>0$; and
$R_n=+\infty$ if and only if $\Lambda=c\delta_1$ for some $c>0$.
\item[(c)] $\cA_n=X_n+\frac\theta2Y_n$ with
$X_n=\sum_{m=3}^{n}[3h_{m-2}+(m-2)(h_{m-3}-h_{m-2})]g_{m-1}$ and
$Y_n=\sum_{m=3}^{n}(m-1)(m-2)(h_{m-3}-h_{m-2})$, and $X_n$, $Y_n$, $\cB_n$ depend on
$\Lambda$ only. If $\Lambda$ is concentrated neither at $0$ nor at $1$, then $Y_n>0$ and
$\cB_n>0$, so $\theta\mapsto R_n$ is affine and strictly increasing on $(0,\infty)$ and
tends to $+\infty$ as $\theta\to\infty$.
\item[(d)] For $n=4$,
\[
R_4=2+\frac{(g_3-3h_2)(g_3+\theta)}{2g_3h_2},\qquad g_3-3h_2=\int x(4-3x)\,\Lambda(dx),
\]
where $g_3=h_0+2h_1$; the integral is the total rate at which, among four lineages, an
event merges at least three of them.
\end{itemize}
\end{thm}

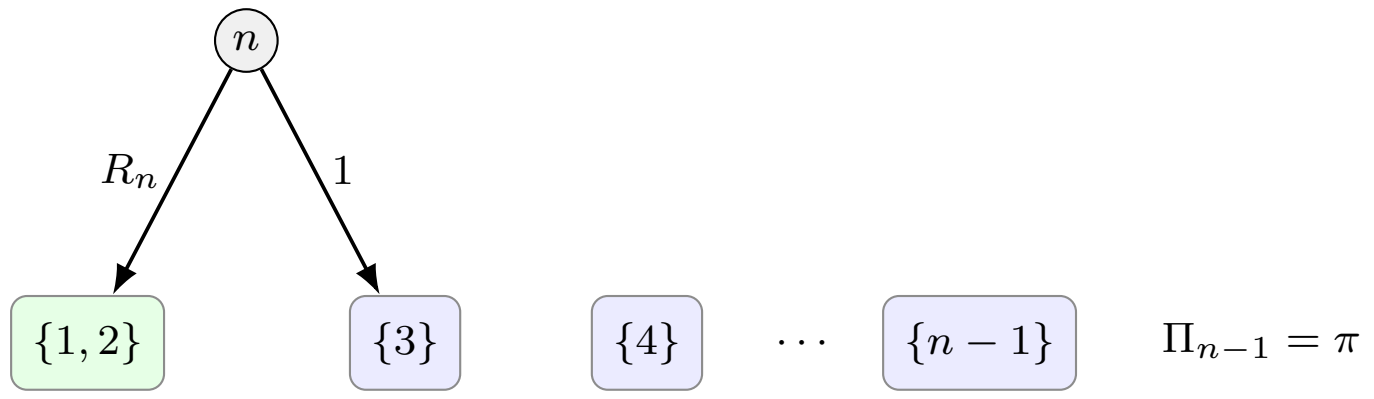
\begin{figure}[htbp]
\centering
\begin{tikzpicture}[font=\footnotesize,
  blk/.style={rounded corners=3pt,draw=black!45,inner sep=4pt,minimum height=0.62cm},
  arr/.style={-{Latex[length=2mm]},line width=0.7pt}]
\node[blk,fill=green!10] (P) at (0,0) {$\{1,2\}$};
\node[blk,fill=blue!8]   (S) at (2.1,0) {$\{3\}$};
\node[blk,fill=blue!8]   (S4) at (3.7,0) {$\{4\}$};
\node at (4.75,0) {$\cdots$};
\node[blk,fill=blue!8]   (Sl) at (5.9,0) {$\{n-1\}$};
\node[draw,circle,inner sep=2pt,fill=black!6] (N) at (1.05,2.0) {$n$};
\draw[arr] (N) -- node[left,pos=0.45]{$R_n$} (P);
\draw[arr] (N) -- node[right,pos=0.45]{$1$} (S);
\node[anchor=west,align=left] at (7.0,0.0) {$\Pi_{n-1}=\pi$};
\end{tikzpicture}
\caption{The edge ratio. Given $\Pi_{n-1}=\pi=\{\{1,2\},\{3\},\dots,\{n-1\}\}$, the
probability that individual $n$ joins the pair $\{1,2\}$ is $R_n$ times the probability that
it joins the singleton $\{3\}$. Under the Ewens sampling formula the ratio is $2$, the ratio
of the block sizes, and under a Gibbs partition of type $\sigma$ it is
$(2-\sigma)/(1-\sigma)$. By Theorem~\ref{thm:edge}(b), $R_n\ge2$ for every $\Lambda$, with
equality only for Kingman's coalescent.}
\label{fig:edge}
\end{figure}

Let $\pi=\{\{1,2\},\{3\},\dots,\{n-1\}\}$ (Figure~\ref{fig:edge}). Adding $n$ to the block $\{1,2\}$ or to the
block $\{3\}$ of $\pi$ gives partitions of $[n]$ with configurations $(3,1^{n-3})$ and
$(2,2,1^{n-4})$, so $R_n$ is the ratio of the conditional probabilities, given
$\Pi_{n-1}=\pi$, that $n$ joins $\{1,2\}$ and that it joins $\{3\}$. For $n\ge4$ we define
the \emph{edge index}
\[
\sigma^{(n)}=\frac{R_n-2}{R_n-1},\qquad\sigma^{(n)}=1\ \text{when}\ R_n=+\infty .
\]
By Theorem~\ref{thm:edge}(b), $\sigma^{(n)}\in[0,1]$, with $\sigma^{(n)}=0$ if and only if
$\Lambda=c\delta_0$ and $\sigma^{(n)}=1$ if and only if $\Lambda=c\delta_1$; when
$R_n<\infty$, $\sigma=\sigma^{(n)}$ is the only solution in $[0,1)$ of
$(2-\sigma)/(1-\sigma)=R_n$.

\begin{proof}[Proof of Theorem~\ref{thm:edge}]
(a) In \eqref{eq:recursion} with $\bn$ an edge configuration, a merger is compatible
with $\bn$ only if it involves at most three lineages, all in one block of $\bn$.
With $\lambda_{n,2}=h_{n-2}$ and $\lambda_{n,3}=h_{n-3}-h_{n-2}$ this gives
\begin{align}
\rho_np_n(1^n)&=nr\,p_{n-1}(1^{n-1}),\label{eq:e1}\\
\rho_np_n(2,1^{n-2})&=(n-2)r\,p_{n-1}(2,1^{n-3})+h_{n-2}\,p_{n-1}(1^{n-1}),\label{eq:e2}\\
\rho_np_n(3,1^{n-3})&=(n-3)r\,p_{n-1}(3,1^{n-4})+3h_{n-2}\,p_{n-1}(2,1^{n-3})\notag\\
&\quad+(h_{n-3}-h_{n-2})\,p_{n-2}(1^{n-2}),\label{eq:e3}\\
\rho_np_n(2,2,1^{n-4})&=(n-4)r\,p_{n-1}(2,2,1^{n-5})+2h_{n-2}\,p_{n-1}(2,1^{n-3}),\label{eq:e4}
\end{align}
where \eqref{eq:e1} and \eqref{eq:e2} hold for $n\ge2$, \eqref{eq:e3} for $n\ge3$ and
\eqref{eq:e4} for $n\ge4$, a term being omitted when its coefficient vanishes. The formula
for $p_n(1^n)$, which is also in the proof of \citet[Lemma~6.3]{Dong07}, follows from \eqref{eq:e1} and
$p_1(1)=1$.
Let $x_n$, $y_n$, $z_n$ be the ratios of $p_n(2,1^{n-2})$, $p_n(3,1^{n-3})$,
$p_n(2,2,1^{n-4})$ to $p_n(1^n)$. Dividing \eqref{eq:e2}--\eqref{eq:e4} by
$p_{n-1}(1^{n-1})$, and using $\rho_np_n(1^n)=nr\,p_{n-1}(1^{n-1})$ and
$p_{n-2}(1^{n-2})/p_{n-1}(1^{n-1})=\rho_{n-1}/((n-1)r)$, both from \eqref{eq:e1}, gives
\begin{align*}
nr\,x_n&=(n-2)r\,x_{n-1}+h_{n-2},\\
nr\,y_n&=(n-3)r\,y_{n-1}+3h_{n-2}\,x_{n-1}+(h_{n-3}-h_{n-2})\rho_{n-1}/((n-1)r),\\
nr\,z_n&=(n-4)r\,z_{n-1}+2h_{n-2}\,x_{n-1},
\end{align*}
with the conventions $x_1=0$, $y_2=0$ and $z_3=0$. Multiplying the first line by
$(n-1)/r$ shows that $\hat x_n=n(n-1)x_n$ satisfies
$\hat x_n=\hat x_{n-1}+(n-1)h_{n-2}/r$ with $\hat x_1=0$, so $\hat x_n=g_n/r$ by
\eqref{eq:g}. Multiplying the second line by $(n-1)(n-2)/r$ and using
$x_{n-1}=g_{n-1}/((n-1)(n-2)r)$ shows that $\hat y_n=n(n-1)(n-2)y_n$ satisfies
$\hat y_n=\hat y_{n-1}+[3h_{n-2}g_{n-1}+(n-2)(h_{n-3}-h_{n-2})\rho_{n-1}]/r^2$ with
$\hat y_2=0$, so $\hat y_n=\cA_n/r^2$. Multiplying the third line by $(n-1)(n-2)(n-3)/r$
shows that $\hat z_n=n(n-1)(n-2)(n-3)z_n$ satisfies
$\hat z_n=\hat z_{n-1}+2(n-3)h_{n-2}g_{n-1}/r^2$ with $\hat z_3=0$, so
$\hat z_n=\cB_n/r^2$. With $r=\theta/2$ this proves the four formulas. The terms of $\cA_n$
are nonnegative, and the term $m=3$, $3h_1g_2+(h_0-h_1)\rho_2$, is positive: if $h_1>0$ its
first summand is, and if $h_1=0$ its second summand equals $h_0\rho_2$.

(b) If $\cB_n=0$ then $h_2g_3=0$; since $g_3\ge h_0>0$, $h_2=0$ and $\Lambda=h_0\delta_1$
by Lemma~\ref{lem:dictionary}. Conversely $\Lambda=c\delta_1$ gives $h_j=0$ for $j\ge1$ and
$\cB_n=0$. Let now $\cB_n>0$. Since $2\cB_n=4\sum_{m=3}^n(m-3)h_{m-2}g_{m-1}$,
\begin{equation}\label{eq:abel}
\begin{aligned}
(n-3)\cA_n-2\cB_n={}&\sum_{m=3}^{n}\nu_m\,h_{m-2}g_{m-1}\\
&+(n-3)\sum_{m=3}^{n}(m-2)(h_{m-3}-h_{m-2})\rho_{m-1},
\end{aligned}
\end{equation}
where $\nu_m=3(n-3)-4(m-3)$. The second sum has nonnegative terms, and its term $m=3$ is
$(h_0-h_1)\rho_2$. For the first sum let $\gamma_m=\binom{m-1}{2}=\sum_{j=0}^{m-3}(j+1)$
and $\varphi_m=h_{m-2}\,g_{m-1}/\gamma_m$, $3\le m\le n$. The quotient $g_{m-1}/\gamma_m$
is the average of $h_0,\dots,h_{m-3}$ with weights $1,\dots,m-2$; passing from $m$ to $m+1$
adds the value $h_{m-2}$, which is at most each of $h_0,\dots,h_{m-3}$, so the average
does not increase. As $h_{m-2}$ does not increase either, $(\varphi_m)$ is nonincreasing.
The weights $\nu_m$ decrease in $m$, so the partial sums $S_M=\sum_{m=3}^{M}\gamma_m\nu_m$
first increase and then decrease, and
\[
S_n=3(n-3)\sum_{m=3}^n\binom{m-1}{2}-12\sum_{m=3}^n\binom{m-1}{3}
=3(n-3)\binom n3-12\binom n4=0,
\]
using $(m-3)\binom{m-1}{2}=3\binom{m-1}{3}$. Hence $S_M\ge0$ for $3\le M\le n$, and by
summation by parts
\[
\sum_{m=3}^{n}\gamma_m\nu_m\varphi_m=\sum_{M=3}^{n-1}S_M(\varphi_M-\varphi_{M+1})+S_n\varphi_n\ge0 .
\]
Therefore $R_n-2=[(n-3)\cA_n-2\cB_n]/\cB_n\ge(n-3)(h_0-h_1)\rho_2/\cB_n\ge0$, and
$R_n>2$ unless $h_0=h_1$, that is, unless $\Lambda=h_0\delta_0$. If $\Lambda=c\delta_0$,
then $h_{m-3}=h_{m-2}$ and $\varphi_m=c^2$ for all $m$, so both sums in \eqref{eq:abel}
vanish and $R_n=2$.

(c) Substitute $\rho_{m-1}=g_{m-1}+(m-1)r$ in \eqref{eq:AB}. The term $m=3$ of $Y_n$ is
$2(h_0-h_1)>0$ when $\Lambda$ is not concentrated at $0$, and $\cB_n>0$ when $\Lambda$ is
not concentrated at $1$ by (b). Then $R_n=(n-3)(X_n+\tfrac\theta2Y_n)/\cB_n$.

(d) For $n=4$, using $g_2=h_0$, $\rho_2=h_0+2r$, $g_3=h_0+2h_1$ and $\rho_3=g_3+3r$,
\begin{align*}
\cA_4&=3h_1h_0+(h_0-h_1)(h_0+2r)+3h_2g_3+2(h_1-h_2)(g_3+3r)\\
&=g_3^2+g_3h_2+2r(g_3-3h_2),
\end{align*}
and $\cB_4=2h_2g_3$, so
$R_4-2=(\cA_4-2\cB_4)/\cB_4=(g_3-3h_2)(g_3+2r)/(2g_3h_2)$. Finally
$1+2(1-x)-3(1-x)^2=x(4-3x)$, and among four lineages the events that merge at least three
of them have total rate $4\lambda_{4,3}+\lambda_{4,4}=\int(4x(1-x)+x^2)\,\Lambda(dx)$.
\end{proof}

The two sums in \eqref{eq:abel} have a direct reading. The second collects the events that
merge three lineages at once, through the rates $\lambda_{m,3}=h_{m-3}-h_{m-2}$; such
events can produce the configuration $(3,1^{n-3})$ but never $(2,2,1^{n-4})$. The first
involves only the rates $\lambda_{m,2}=h_{m-2}$ at which two given lineages merge, and it is
nonnegative because these rates do not increase with the number $m$ of lineages.

\begin{cor}\label{cor:nonpositive}
Let $n\ge4$ and $\sigma\in[-\infty,1)$. If $\Pi_n$ is of Gibbs type $\sigma$, then
$R_n=(2-\sigma)/(1-\sigma)$ if $\sigma>-\infty$, and $R_n=1$ if $\sigma=-\infty$.
Consequently:
\begin{itemize}
\item[(a)] if $\Pi_n$ is of Gibbs type $\sigma\in[-\infty,0]$, then $\Lambda=c\delta_0$
for some $c>0$, and $\sigma=0$;
\item[(b)] if $\Lambda=c\delta_1$, then $\Pi_n$ is of Gibbs type $\sigma$ for no
$\sigma\in[-\infty,1)$.
\end{itemize}
\end{cor}

\begin{proof}
The configurations $(3,1^{n-3})$ and $(2,2,1^{n-4})$ both have $n-2$ blocks, and
$w_\sigma(3,1^{n-3})=(1-\sigma)(2-\sigma)$, $w_\sigma(2,2,1^{n-4})=(1-\sigma)^2$ (both
equal to $1$ if $\sigma=-\infty$). With $V=V_{n,n-2}$, $p_n(3,1^{n-3})=V(1-\sigma)(2-\sigma)$
and $p_n(2,2,1^{n-4})=V(1-\sigma)^2$; $V>0$ because $p_n(3,1^{n-3})>0$ by
Theorem~\ref{thm:edge}(a). This gives the value of $R_n$, which is finite. For
$\sigma\in(-\infty,0]$, $(2-\sigma)/(1-\sigma)=1+1/(1-\sigma)\le2$, with equality only at
$\sigma=0$, and for $\sigma=-\infty$, $R_n=1$. Theorem~\ref{thm:edge}(b) leaves only $R_n=2$,
$\sigma=0$ and $\Lambda=c\delta_0$. If $\Lambda=c\delta_1$, then $R_n=+\infty$, which no
type allows.
\end{proof}

Part (d) of Theorem~\ref{thm:edge} yields a statistic, computed from the sample alone, whose
mean vanishes for every $\theta$ under Kingman's coalescent and is positive for every
$\theta$ under every other $\Lambda$-coalescent.

\begin{cor}\label{cor:ustat}
Let $n\ge4$, and let $U_{31}$ and $U_{22}$ be the proportions of the $\binom n4$ subsets of
$[n]$ with four elements on which $\Pi_n$ induces a partition with configuration $(3,1)$
and $(2,2)$, respectively. Then
\[
\E\bigl[3U_{31}-8U_{22}\bigr]=\frac{6\theta\,(g_3+\theta)}{\rho_2\rho_3\rho_4}
\int x(4-3x)\,\Lambda(dx).
\]
The right-hand side is the same for all $n\ge4$; it vanishes for every $\theta$ if
$\Lambda=c\delta_0$, and it is positive for every $\theta$ otherwise.
\end{cor}

\begin{proof}
By exchangeability and consistency, the restriction of $\Pi_n$ to a subset with four
elements has the law of $\Pi_4$. There are $4$ partitions of $[4]$ with configuration
$(3,1)$ and $3$ with configuration $(2,2)$, so $\E[U_{31}]=4p_4(3,1)$,
$\E[U_{22}]=3p_4(2,2)$ and $\E[3U_{31}-8U_{22}]=12(p_4(3,1)-2p_4(2,2))
=p_4(1^4)(\cA_4-2\cB_4)/(2r^2)$ by Theorem~\ref{thm:edge}(a). By the proof of
Theorem~\ref{thm:edge}(d), $\cA_4-2\cB_4=(g_3-3h_2)(g_3+2r)$, and
$p_4(1^4)=24r^3/(\rho_2\rho_3\rho_4)$. Since $x(4-3x)>0$ on $(0,1]$, the integral is
positive unless $\Lambda=c\delta_0$.
\end{proof}

\paragraph{The Bolthausen--Sznitman coalescent}
For $\Lambda$ the Lebesgue measure on $[0,1]$ \citep{Bolthausen98}, $h_j=1/(j+1)$, so
$g_m=m-1$ by \eqref{eq:g} and $h_{m-3}-h_{m-2}=1/((m-2)(m-1))$. The terms of $\cA_n$ in
\eqref{eq:AB} are $4(m-2)/(m-1)+\theta/2$ and those of $\cB_n$ are $2(m-4)+4/(m-1)$, so
that, with the harmonic numbers $H_m=\sum_{i=1}^m1/i$,
\begin{equation}\label{eq:BS}
\cA_n=(n-2)\Bigl(4+\frac\theta2\Bigr)-4H_{n-1}+4,\qquad\cB_n=(n-3)(n-4)+4H_{n-1}-6 .
\end{equation}
In particular $R_4=(14+3\theta)/4$, $R_5=(46+9\theta)/13$ and $R_6=3(163+30\theta)/137$.

\section{Sufficiency of the number of alleles}\label{sec:sufficiency}

Fix $\Lambda$ and $n$, and consider the laws of $\Pi_n$ for $\theta\in(0,\infty)$. We say
that $K_n$ is \emph{sufficient for $\theta$} if, for each $k$, the conditional law of
$\Pi_n$ given $K_n=k$ is the same for all $\theta$ with $\PP_\theta(K_n=k)>0$.

\begin{thm}\label{thm:sufficiency}
Fix a finite nonzero measure $\Lambda$ on $[0,1]$.
\begin{itemize}
\item[(i)] If $n\le3$, $K_n$ is sufficient for $\theta$.
\item[(ii)] If $n\ge4$, $K_n$ is sufficient for $\theta$ if and only if $\Lambda=c\delta_0$
or $\Lambda=c\delta_1$ for some $c>0$.
\item[(iii)] Let $n\ge4$ and $\theta_1\ne\theta_2$. If $\Pi_n$ is of Gibbs type $\sigma$
both under $\theta_1$ and under $\theta_2$, for the same $\sigma\in[-\infty,1)$, then
$\Lambda=c\delta_0$ for some $c>0$, and $\sigma=0$.
\item[(iv)] For $n\ge4$ and every $\theta>0$,
\begin{align*}
\PP\bigl(\Pi_n\text{ has configuration }(3,1^{n-3})\bigm|K_n=n-2\bigr)
&=\frac{R_n}{R_n+3(n-3)/4}\\
&\ge\frac{8}{3n-1},
\end{align*}
with equality if and only if $\Lambda=c\delta_0$. The probability equals $1$ if and only if
$\Lambda=c\delta_1$; if $\Lambda$ is concentrated neither at $0$ nor at $1$, it is a
strictly increasing function of $\theta$.
\end{itemize}
\end{thm}

\begin{proof}
(i) For $n\le3$ the configuration is determined by $K_n$, and given the configuration the
partition is uniform among the partitions with that configuration, by exchangeability.

(iv) The partitions of $[n]$ with $n-2$ blocks have configuration $(3,1^{n-3})$, and there
are $\binom n3$ of them, or $(2,2,1^{n-4})$, and there are $3\binom n4$ of them. The
probability is therefore
$\binom n3p_n(3,1^{n-3})/[\binom n3p_n(3,1^{n-3})+3\binom n4p_n(2,2,1^{n-4})]
=R_n/(R_n+3(n-3)/4)$, read as $1$ when $R_n=+\infty$; $\PP(K_n=n-2)>0$ since
$p_n(3,1^{n-3})>0$. The map $x\mapsto x/(x+3(n-3)/4)$ is increasing, and its value at $x=2$
is $8/(3n-1)$. The assertions follow from Theorem~\ref{thm:edge}(b) and (c).

(ii) If $\Lambda=c\delta_0$, the conditional law of $\Pi_n$ given $K_n=k$ is proportional to
$\prod_j(n_j-1)!$ by Corollary~\ref{cor:ewens}, whatever $\theta$. If
$\Lambda=c\delta_1$, then $\lambda_{b,m}=c\,\mathbf 1\{m=b\}$, and \eqref{eq:recursion}
shows by induction on $n$ that $p_n(\bn)=0$ unless the configuration of $\bn$ is
$(n-k+1,1^{k-1})$ for some $k$: removing a singleton from a configuration with two parts of
size at least $2$ leaves two such parts, and the only merger term is that of $m=n$, which
requires $\bn=(n)$. Given $K_n=k$ the configuration is thus $(n-k+1,1^{k-1})$, and the
partition is uniform among those with this configuration, whatever $\theta$. If $\Lambda$
is concentrated neither at $0$ nor at $1$, the conditional probability in (iv) is strictly
increasing in $\theta$, so $K_n$ is not sufficient.

(iii) If $\sigma\le0$, Corollary~\ref{cor:nonpositive}(a) applies. If $\sigma\in(0,1)$,
Corollary~\ref{cor:nonpositive} gives $R_n=(2-\sigma)/(1-\sigma)\in(2,\infty)$ under both
$\theta_1$ and $\theta_2$; then $\Lambda$ is concentrated neither at $0$ nor at $1$, and
Theorem~\ref{thm:edge}(c) gives $R_n(\theta_1)\ne R_n(\theta_2)$, a contradiction.
\end{proof}

Theorem~\ref{thm:sufficiency}(iii) is the form of the rigidity in which the Gibbs type is
the same for two values of $\theta$, as in the Ewens family, where it is $0$ for every
$\theta$; it holds at every $n\ge4$. Section~\ref{sec:rigidity} treats a single value of $\theta$.

\section{Gibbs rigidity at a fixed \texorpdfstring{$\theta$}{theta}}\label{sec:rigidity}

\begin{thm}\label{thm:rigidity}
Let $\Lambda$ be a finite nonzero measure on $[0,1]$ and let $\theta>0$ be fixed. The
following are equivalent.
\begin{itemize}
\item[(i)] $\Pi_6$ is of Gibbs type $\sigma$ for some $\sigma\in[-\infty,1)$.
\item[(ii)] $\Lambda=c\delta_0$ for some $c>0$.
\item[(iii)] For every $n\ge1$, $\Pi_n$ follows the Ewens sampling formula with parameter
$\theta/c$; in particular $\Pi_n$ is of Gibbs type $0$ for every $n$.
\end{itemize}
\end{thm}

\begin{cor}\label{cor:rigidity}
Let $\theta>0$ and $n\ge6$. Then $\Pi_n$ is of Gibbs type $\sigma$ for some
$\sigma\in[-\infty,1)$ if and only if $\Lambda=c\delta_0$. The allelic partition $\Pi$ of
$\N$ is of Gibbs type if and only if $\Lambda=c\delta_0$.
\end{cor}

\begin{proof}
By Lemma~\ref{lem:addition}, a partition of $[n]$ of Gibbs type $\sigma$, $n\ge6$, restricts
to a partition of $[6]$ of the same type; apply Theorem~\ref{thm:rigidity}.
\end{proof}

\subsection{Sample size five}\label{sec:five}

\begin{lem}\label{lem:five}
Let $(\Pi_n)_{n\ge1}$ be the restrictions of an exchangeable random partition of $\N$ and
let $\sigma\in[-\infty,1)$. Then $\Pi_5$ is of Gibbs type $\sigma$ if and only if $\Pi_4$ is
of Gibbs type $\sigma$ and $u_5(3,1,1)=u_5(2,2,1)$.
\end{lem}

\begin{proof}
The condition is necessary by Lemma~\ref{lem:addition} and the definition. Conversely, the
configurations of $[5]$ with $k=1,4,5$ blocks are unique, those with $k=3$ are $(3,1,1)$ and
$(2,2,1)$, and those with $k=2$ are $(4,1)$ and $(3,2)$, so it remains to show
$u_5(4,1)=u_5(3,2)$. Put $V=V_{4,2}$ and $W=u_5(3,1,1)=u_5(2,2,1)$, and let $\sigma>-\infty$.
By consistency,
\[
p_4(3,1)=p_5(4,1)+p_5(3,2)+p_5(3,1,1),\qquad p_4(2,2)=2p_5(3,2)+p_5(2,2,1),
\]
where $p_4(3,1)=V(1-\sigma)(2-\sigma)$, $p_4(2,2)=V(1-\sigma)^2$,
$p_5(3,1,1)=W(1-\sigma)(2-\sigma)$ and $p_5(2,2,1)=W(1-\sigma)^2$. The second identity gives
$p_5(3,2)=(V-W)(1-\sigma)^2/2$, and then the first gives
$p_5(4,1)=(V-W)(1-\sigma)(3-\sigma)/2$. Dividing by $w_\sigma(3,2)=(1-\sigma)^2(2-\sigma)$
and by $w_\sigma(4,1)=(1-\sigma)(2-\sigma)(3-\sigma)$,
$u_5(3,2)=u_5(4,1)=(V-W)/(2(2-\sigma))$. For $\sigma=-\infty$ all weights are $1$, and the
same identities give $p_5(3,2)=p_5(4,1)=(V-W)/2$.
\end{proof}

For the allelic partition with $\Lambda\ne c\delta_1$ and $\sigma>-\infty$, the paragraph
after Definition~\ref{def:gibbs} and Lemma~\ref{lem:five} show that $\Pi_5$ is of Gibbs type
$\sigma$ if and only if $R_4=R_5=(2-\sigma)/(1-\sigma)$.

\begin{prop}\label{prop:sharp}
\begin{itemize}
\item[(a)] If $\Lambda$ is concentrated neither at $0$ nor at $1$, then for every
$\theta>0$ the partition $\Pi_4$ is of Gibbs type $\sigma^{(4)}\in(0,1)$ and of no other
type.
\item[(b)] Let $\Lambda$ be the Lebesgue measure on $[0,1]$. Then
$R_4-R_5=(3\theta-2)/52$, so $\Pi_5$ is of Gibbs type for some $\sigma$ if and only if
$\theta=2/3$. For $\theta=2/3$, $R_4=R_5=4$ and $R_6=549/137$: $\Pi_4$ and $\Pi_5$ are of
Gibbs type $2/3$, and $\Pi_6$ is of Gibbs type for no $\sigma\in[-\infty,1)$.
\end{itemize}
\end{prop}

\begin{proof}
(a) By the paragraph after Definition~\ref{def:gibbs}, only the configurations $(3,1)$ and
$(2,2)$ constrain the type of $\Pi_4$. By Theorem~\ref{thm:edge}(b), $R_4\in(2,\infty)$, so
$(1-\sigma)p_4(3,1)=(2-\sigma)p_4(2,2)$ holds exactly for $\sigma=\sigma^{(4)}$, while
$\sigma=-\infty$ would need $R_4=1$. (b) By \eqref{eq:BS},
$R_4-R_5=(14+3\theta)/4-(46+9\theta)/13=(3\theta-2)/52$. In $[4]$ the only possible type is
$\sigma^{(4)}$, so by Lemma~\ref{lem:five} $\Pi_5$ is of Gibbs type exactly when $R_4=R_5$.
For $\theta=2/3$, \eqref{eq:BS} gives $R_4=R_5=4=(2-\sigma)/(1-\sigma)$ with $\sigma=2/3$,
and $R_6=3(163+20)/137=549/137\ne4$. By Lemma~\ref{lem:addition} a Gibbs type of $\Pi_6$
would be a Gibbs type of $\Pi_4$, hence $2/3$, which needs $R_6=4$.
\end{proof}

Since replacing $(\Lambda,\theta)$ by $(c\Lambda,c\theta)$ does not change the law of
$\Pi_n$, part (b) applies to $\Lambda=(3\theta/2)\,\mathrm{Leb}$ for every $\theta>0$, where
$\mathrm{Leb}$ is the Lebesgue measure on $[0,1]$: for this measure $\Pi_5$ is of Gibbs type
$2/3$ and $\Pi_6$ is of Gibbs type for no $\sigma$. The sample size six in
Theorem~\ref{thm:rigidity} is therefore minimal for each fixed $\theta$, and by part (a)
four does not suffice for any $\Lambda$ concentrated neither at $0$ nor at $1$.

\subsection{Proof of Theorem~\ref{thm:rigidity}}\label{sec:proof-rigidity}

Under Gibbs type $\sigma>-\infty$ on $[n]$, Corollary~\ref{cor:nonpositive} gives the
\emph{edge identity}
\begin{equation}\label{eq:edge-identity}
(1-\sigma)\,p_n(3,1^{n-3})=(2-\sigma)\,p_n(2,2,1^{n-4}).
\end{equation}

\begin{lem}\label{lem:affine}
Let $n\ge4$ and $\sigma\in(0,1)$. By Theorem~\ref{thm:edge}(a), \eqref{eq:edge-identity}
is equivalent to $\Delta_n=0$, where $\Delta_n=(1-\sigma)(n-3)\cA_n-(2-\sigma)\cB_n$. With
$\theta,h_0,\dots,h_{n-3}$ fixed, $\Delta_n$ is an affine function of $h_{n-2}$ with slope
\[
\kappa_n=-(n-3)\bigl[(1+\sigma)g_{n-1}+(1-\sigma)(n-2)\rho_{n-1}\bigr]<0 .
\]
In particular $\Delta_4=0$ if and only if
\begin{equation}\label{eq:H2}
h_2=\frac{(1-\sigma)\,g_3\,(g_3+\theta)}{(3-\sigma)g_3+3(1-\sigma)\theta},\qquad g_3=h_0+2h_1 .
\end{equation}
\end{lem}

\begin{proof}
The moment $h_{n-2}$ enters $\cA_n$ and $\cB_n$ only through their terms $m=n$, since
$g_{n-1}$ and $\rho_{n-1}$ involve $h_0,\dots,h_{n-3}$ only. Hence
$\partial\cA_n/\partial h_{n-2}=3g_{n-1}-(n-2)\rho_{n-1}$ and
$\partial\cB_n/\partial h_{n-2}=2(n-3)g_{n-1}$, and
$\kappa_n=(n-3)[3(1-\sigma)g_{n-1}-2(2-\sigma)g_{n-1}-(1-\sigma)(n-2)\rho_{n-1}]$, which is
the stated value; it is negative because $g_{n-1}\ge h_0>0$ and $\rho_{n-1}>0$. For $n=4$,
the proof of Theorem~\ref{thm:edge}(d) gives
$\Delta_4=(1-\sigma)\cA_4-2(2-\sigma)g_3h_2=(1-\sigma)g_3(g_3+2r)-h_2[(3-\sigma)g_3+6(1-\sigma)r]$,
and $r=\theta/2$.
\end{proof}

In the rest of this subsection and in Appendix~\ref{app:certificate} we take $h_0=1$; the
proof of Theorem~\ref{thm:rigidity} reduces to this case by the scaling of
Section~\ref{sec:model}. Let $\sigma\in(0,1)$, $\theta>0$ and $h_1\in(0,1)$. By
Lemma~\ref{lem:affine}, $\Delta_4=0$ has the unique solution $h_2=\eta_2(\sigma,\theta,h_1)$
given by \eqref{eq:H2}. With $h_2=\eta_2$, $\Delta_5$ is affine in $h_3$ with slope
$\kappa_5$, and $\Delta_5=0$ defines $h_3=\eta_3(\sigma,\theta,h_1)$; with $h_2=\eta_2$ and
$h_3=\eta_3$, $\Delta_6=0$ defines $h_4=\eta_4(\sigma,\theta,h_1)$. The slopes $\kappa_5$
and $\kappa_6$ are here evaluated at $h_2=\eta_2$ and $h_3=\eta_3$; their explicit forms,
obtained in the computation behind Lemma~\ref{lem:certificate} and listed in
Appendix~\ref{app:certificate}, are negative for all $\sigma\in(0,1)$, $\theta>0$ and
$h_1\in(0,1)$, so that $\eta_2$, $\eta_3$, $\eta_4$ are rational functions defined there.
Let $\Def_6(\sigma,\theta,h_1)$ be the rational function $u_6(3,2,1)-u_6(2,2,2)$ obtained
from \eqref{eq:recursion}, \eqref{eq:dictionary} and \eqref{eq:g} with the moments
$(1,h_1,\eta_2,\eta_3,\eta_4)$. Both configurations have three blocks, so a partition of
$[6]$ of Gibbs type $\sigma$ has $u_6(3,2,1)=u_6(2,2,2)$. The letters $a$, $b$, $c$ denote,
in Lemma~\ref{lem:certificate} and Appendix~\ref{app:certificate}, the coordinates
\[
a=\frac{\sigma}{1-\sigma},\qquad b=\theta,\qquad c=\frac{h_1}{1-h_1},
\]
which map $(0,1)\times(0,\infty)\times(0,1)$ bijectively onto $(0,\infty)^3$, with
inverse $\sigma=a/(1+a)$, $\theta=b$, $h_1=c/(1+c)$.

\begin{lem}\label{lem:certificate}
With $h_0=1$ and the coordinates $a,b,c$ above,
\begin{equation}\label{eq:certificate}
\Def_6(\sigma,\theta,h_1)=\frac{a^2b^2(1+a)^3(1+c)(1+3c)^4(1+b+3c+bc)\,Q}
{(2+a)(1+b)(2+3b+6c+3bc)\,F_4F_8F_{16}},
\end{equation}
where the polynomials $Q,F_4,F_8,F_{16}\in\mathbb Z[a,b,c]$ have total degrees
$18$, $4$, $8$, $16$, have $530$, $15$, $70$, $396$ terms, and have only positive
coefficients. In particular
$\Def_6(\sigma,\theta,h_1)>0$ for all $\sigma\in(0,1)$, $\theta>0$ and $h_1\in(0,1)$.
\end{lem}

\begin{proof}
The left-hand side of \eqref{eq:certificate} is obtained from $\sigma$, $\theta$ and $h_1$
by the finitely many field operations described above: $\eta_2$, $\eta_3$ and $\eta_4$ from
Lemma~\ref{lem:affine}, then $p_m(\bn)$ for every composition $\bn$ of every $m\le6$ from
\eqref{eq:recursion} with the moments $(1,h_1,\eta_2,\eta_3,\eta_4)$, and finally
$p_6(3,2,1)/w_\sigma(3,2,1)-p_6(2,2,2)/w_\sigma(2,2,2)$. Carrying out these operations by
exact arithmetic in the field $\mathbb Q(\sigma,\theta,h_1)$ of rational functions,
substituting $\sigma=a/(1+a)$, $\theta=b$, $h_1=c/(1+c)$ and reducing to lowest terms gives
\eqref{eq:certificate}. The computation was carried out twice, once with the edge
identities solved through the closed forms of Theorem~\ref{thm:edge} and once directly
from \eqref{eq:recursion}, with the same result. The coefficients of $Q$, $F_4$, $F_8$ and
$F_{16}$, listed in the supplementary file \texttt{certificate-n6.json}, are positive
integers, so every factor on the right of \eqref{eq:certificate} is positive on
$(0,\infty)^3$.
\end{proof}

\begin{proof}[Proof of Theorem~\ref{thm:rigidity}]
(ii)$\Rightarrow$(iii) is Corollary~\ref{cor:ewens}, and (iii)$\Rightarrow$(i) is
immediate. Assume (i). If $\sigma\le0$, Corollary~\ref{cor:nonpositive}(a) gives (ii).
Suppose $\sigma\in(0,1)$. By Corollary~\ref{cor:nonpositive}, $\Lambda$ is not a multiple of
$\delta_1$, and it is not a multiple of $\delta_0$, because then
$R_6=2<(2-\sigma)/(1-\sigma)$. Hence $0<h_1<h_0$ by Lemma~\ref{lem:dictionary}. Replacing $(\Lambda,\theta)$ by
$(\Lambda/h_0,\theta/h_0)$ does not change the law of $\Pi_6$ and divides every $h_j$ by
$h_0$; it therefore replaces $h_0$ by $1$, $\theta$ by $\theta/h_0$ and $h_1$ by
$h_1/h_0\in(0,1)$, and in the coordinates of Lemma~\ref{lem:certificate} the point is
$a=\sigma/(1-\sigma)$, $b=\theta/h_0$, $c=h_1/(h_0-h_1)$. From now on $\theta$ and $h_j$
denote these new values, so that $h_0=1$. By Lemma~\ref{lem:addition}, $\Pi_4$, $\Pi_5$
and $\Pi_6$ are of Gibbs type $\sigma$, so $\Delta_4=\Delta_5=\Delta_6=0$ at the moments of
$\Lambda$. By Lemma~\ref{lem:affine}, whose slopes are now computed at the moments of
$\Lambda$ and are negative, the moments $h_2,h_3,h_4$ of $\Lambda$ are the values of
$\eta_2,\eta_3,\eta_4$ at $(\sigma,\theta,h_1)$. The rational function $\Def_6$ is built
from these by field operations whose denominators, evaluated at this point, are the rates
$\rho_2,\dots,\rho_6$ of $\Lambda$, the slopes $\kappa_4,\kappa_5,\kappa_6$ at the moments
of $\Lambda$, and powers of $1-\sigma$, $2-\sigma$ and $1-h_1$, none of which vanishes; the
denominator of \eqref{eq:certificate} is positive there as well. Both sides of
\eqref{eq:certificate} may therefore be evaluated at the point, and
$u_6(3,2,1)-u_6(2,2,2)=\Def_6(\sigma,\theta,h_1)>0$ by Lemma~\ref{lem:certificate}. This
contradicts the Gibbs identity $u_6(3,2,1)=u_6(2,2,2)$.
\end{proof}

\section{The edge index for large samples}\label{sec:limits}

\begin{lem}\label{lem:rare}
For every $\Lambda$ and every $\theta>0$, $\PP(K_n\ge n-2)$ is nonincreasing in $n$ and
tends to $0$ as $n\to\infty$.
\end{lem}

\begin{proof}
Let $\Pi$ be the allelic partition of $\N$, whose restriction to $[n]$ is $\Pi_n$. Passing
from $[n]$ to $[n+1]$, the element $n+1$ joins a block or forms a new one, so $n-K_n$ is
nondecreasing in $n$; it converges to $\sum_B(|B|-1)\in[0,\infty]$, the sum over the
blocks $B$ of $\Pi$. Hence the events $\{K_n\ge n-2\}=\{n-K_n\le2\}$ decrease, and
$\PP(K_n\ge n-2)$ decreases to $\PP(E)$, where $E=\{\sum_B(|B|-1)\le2\}$. Write $i\sim j$
when $i$ and $j$ lie in the same block of $\Pi$. A partition of $[N]$ with at least $N-2$
blocks has at most three pairs $\{i,j\}$ of elements in a common block, so
$\sum_{1\le i<j\le N}\PP(i\sim j,\,E)\le3$. The event $E$ is invariant under finite
permutations of $\N$ and $\Pi$ is exchangeable, so every term of this sum equals
$\PP(1\sim2,\,E)$, and $\binom N2\PP(1\sim2,\,E)\le3$ for every $N$. Thus
$\PP(i\sim j,\,E)=0$ for all $i<j$, and up to a null set $E$ is contained in the event that
every block of $\Pi$ is a singleton, whose probability is $\lim_np_n(1^n)$. By
Theorem~\ref{thm:edge}(a), $p_n(1^n)=\prod_{m=2}^n\bigl(1+2g_m/(m\theta)\bigr)^{-1}$, and
$g_m\ge h_0>0$, so $\sum_m2g_m/(m\theta)=\infty$ and $p_n(1^n)\to0$.
\end{proof}

By Lemma~\ref{lem:rare}, the edge index describes the law of $\Pi_n$ on events whose
probability tends to $0$, not the typical configurations of large samples.
Proposition~\ref{prop:limits} gives its limits, and the comparison that follows it, with the
block counts of \citet{BBS2007}, contrasts the two regimes.

\begin{prop}\label{prop:limits}
Fix $\theta>0$ and let $n\to\infty$.
\begin{itemize}
\item[(a)] If $\Lambda(\{0\})>0$, then $R_n\to2$ and $\sigma^{(n)}\to0$.
\item[(b)] Let $1<\alpha<2$, and let $\Lambda$ have a density $f$ on $[0,1]$ with
$f(x)\sim Cx^{1-\alpha}$ as $x\downarrow0$, for some $C>0$; for instance
$\Lambda=\Beta(2-\alpha,\alpha)$, the probability measure with density
$x^{1-\alpha}(1-x)^{\alpha-1}/\mathrm B(2-\alpha,\alpha)$. Then
$R_n\to\alpha(5-\alpha)/(2\alpha-1)$ and
\[
\sigma^{(n)}\to\sigma_\star(\alpha)=\frac{(2-\alpha)(1+\alpha)}{1+3\alpha-\alpha^2};
\]
the limits depend neither on $\theta$ nor on $C$.
\item[(c)] If $\Lambda$ is the Lebesgue measure on $[0,1]$, then for every $n\ge4$
\[
R_n=\frac{(n-3)\bigl[(n-2)(4+\theta/2)-4H_{n-1}+4\bigr]}{(n-3)(n-4)+4H_{n-1}-6},
\]
so that $R_n=4+\theta/2+O(\log n/n)$ and $\sigma^{(n)}\to(4+\theta)/(6+\theta)$.
\item[(d)] Let $0<\alpha<1$, and let $\Lambda$ have a density $f$ on $[0,1]$ with
$f(x)\sim Cx^{1-\alpha}$ as $x\downarrow0$, for some $C>0$. Then
$R_n\sim\theta\alpha(2-\alpha)\,n^{1-\alpha}/(2C\,\Gamma(2-\alpha))$ and $\sigma^{(n)}\to1$.
For $\Lambda=\Beta(2-\alpha,\alpha)$, $R_n\sim\tfrac\theta2(2-\alpha)\Gamma(\alpha+1)\,n^{1-\alpha}$.
\end{itemize}
\end{prop}

\begin{proof}
We use two elementary facts. If $v_m\sim\ell m^\beta$ with $\ell>0$ and $\beta>-1$, then
$\sum_{m\le n}v_m\sim\ell n^{\beta+1}/(\beta+1)$, by comparison with
$\int_0^nt^\beta\,dt$; and if $v_m=O(m^\beta)$, then $\sum_{m\le n}v_m=O(n^{\beta+1})$ for
$\beta>-1$, $O(\log n)$ for $\beta=-1$ and $O(1)$ for $\beta<-1$.

(a) Let $c_0=\Lambda(\{0\})$, so that $h_j\to c_0$, and let $\gamma_m=\binom{m-1}2$. The
sequence $(\gamma_m)$ increases to $\infty$, and by \eqref{eq:g}
\[
\frac{g_m-g_{m-1}}{\gamma_{m+1}-\gamma_m}=\frac{(m-1)h_{m-2}}{m-1}=h_{m-2}\to c_0,
\]
so $g_{m-1}/\gamma_m\to c_0$ by the Stolz--Ces\`aro theorem. Hence
$h_{m-2}g_{m-1}\sim c_0^2m^2/2$, $\sum_{m=3}^n3h_{m-2}g_{m-1}\sim c_0^2n^3/2$ and
$\cB_n\sim c_0^2n^4/4$. The remaining part of $\cA_n$ is $o(n^3)$: since
$g_{m-1}\le h_0\gamma_m$, we have $(m-2)\rho_{m-1}\le(h_0/2+r)m^3$, and
$\sum_m(h_{m-3}-h_{m-2})\le h_0$, so for $\varepsilon\in(0,1)$
\[
\sum_{m=3}^{n}(m-2)(h_{m-3}-h_{m-2})\rho_{m-1}
\le(h_0/2+r)\bigl[\varepsilon^3n^3h_0+n^3(h_{\lfloor\varepsilon n\rfloor-2}-c_0)\bigr],
\]
which is at most $2(h_0/2+r)h_0\varepsilon^3n^3$ for large $n$. Therefore
$\cA_n\sim c_0^2n^3/2$ and $R_n=(n-3)\cA_n/\cB_n\to2$.

For (b) and (d), let $0<\alpha<2$ and let $f$ be as in (b). We first show that, for
$k\in\{0,1\}$,
\begin{equation}\label{eq:abelian}
\int_{[0,1]}x^k(1-x)^j\,\Lambda(dx)\sim C\,\Gamma(k+2-\alpha)\,j^{\alpha-2-k}\qquad(j\to\infty).
\end{equation}
Choose $\delta\in(0,1)$ with $f(x)\le2Cx^{1-\alpha}$ for $0<x<\delta$. The integral over
$[\delta,1]$ is at most $(1-\delta)^j\Lambda([\delta,1])$, which decays exponentially in
$j$. On $(0,\delta)$ the substitution $x=t/j$ gives
\[
j^{k+2-\alpha}\int_0^\delta x^k(1-x)^jf(x)\,dx
=\int_0^{\delta j}t^{k+1-\alpha}\Bigl(1-\frac tj\Bigr)^j\phi\Bigl(\frac tj\Bigr)dt,
\]
where $\phi(x)=x^{\alpha-1}f(x)\to C$ as $x\downarrow0$ and $\phi\le2C$ on $(0,\delta)$.
Since $0\le(1-t/j)^j\le e^{-t}$ for $0\le t\le j$ and $(1-t/j)^j\to e^{-t}$, dominated
convergence, with the integrable bound $2Ct^{k+1-\alpha}e^{-t}$ ($k+1-\alpha>-1$), gives
the limit $C\int_0^\infty t^{k+1-\alpha}e^{-t}\,dt=C\,\Gamma(k+2-\alpha)$. With $k=0$ and
$k=1$, and $L=C\,\Gamma(2-\alpha)$, \eqref{eq:abelian} gives $h_j\sim Lj^{\alpha-2}$ and
$h_j-h_{j+1}=\int x(1-x)^j\Lambda(dx)\sim(2-\alpha)Lj^{\alpha-3}$, because
$\Gamma(3-\alpha)=(2-\alpha)\Gamma(2-\alpha)$. Hence $h_{m-2}\sim Lm^{\alpha-2}$,
$h_{m-3}-h_{m-2}\sim(2-\alpha)Lm^{\alpha-3}$ and, by \eqref{eq:g} and the first elementary
fact, $g_{m-1}\sim Lm^{\alpha}/\alpha$. The terms of $\cB_n$ are then asymptotic to
$2L^2m^{2\alpha-1}/\alpha$, and $\cB_n\sim L^2n^{2\alpha}/\alpha^2$.

(b) For $1<\alpha<2$, $\rho_{m-1}\sim g_{m-1}$, so the terms of $\cA_n$ are asymptotic to
$3L^2m^{2\alpha-2}/\alpha+(2-\alpha)L^2m^{2\alpha-2}/\alpha$, and
$\cA_n\sim(5-\alpha)L^2n^{2\alpha-1}/(\alpha(2\alpha-1))$. Hence
$R_n=(n-3)\cA_n/\cB_n\to R=\alpha(5-\alpha)/(2\alpha-1)$. Then
$R-2=(2-\alpha)(1+\alpha)/(2\alpha-1)$ and $R-1=(1+3\alpha-\alpha^2)/(2\alpha-1)$, which
gives $\sigma_\star(\alpha)$. For $\Beta(2-\alpha,\alpha)$ the density is
$x^{1-\alpha}(1-x)^{\alpha-1}/\mathrm B(2-\alpha,\alpha)$, so (b) applies.

(c) This is \eqref{eq:BS}. The difference between the numerator and $(4+r)$ times the
denominator is
$2(4+r)(n-3)-(4+r)(4H_{n-1}-6)-(n-3)(4H_{n-1}-4)=O(n\log n)$, while the denominator is
asymptotic to $n^2$; and $(2+r)/(3+r)=(4+\theta)/(6+\theta)$.

(d) For $0<\alpha<1$, $g_{m-1}=o(m)$, so $\rho_{m-1}\sim mr$, and the terms
$(m-2)(h_{m-3}-h_{m-2})\rho_{m-1}$ are asymptotic to $r(2-\alpha)Lm^{\alpha-1}$, with sum
asymptotic to $r(2-\alpha)Ln^\alpha/\alpha$. The terms $3h_{m-2}g_{m-1}$ are
$O(m^{2\alpha-2})$, with sum $O(n^{2\alpha-1}+\log n)=o(n^\alpha)$ because $2\alpha-1<\alpha$.
Hence $\cA_n\sim r(2-\alpha)Ln^\alpha/\alpha$ and
$R_n=(n-3)\cA_n/\cB_n\sim r\alpha(2-\alpha)n^{1-\alpha}/L\to\infty$, which is the stated
equivalent with $r=\theta/2$. For $\Beta(2-\alpha,\alpha)$,
$C=1/\mathrm B(2-\alpha,\alpha)=1/(\Gamma(2-\alpha)\Gamma(\alpha))$, so $L=1/\Gamma(\alpha)$
and $r\alpha(2-\alpha)/L=r(2-\alpha)\Gamma(\alpha+1)$.
\end{proof}

\paragraph{Comparison with the block counts}
For $\Lambda=\Beta(2-\alpha,\alpha)$ with $1<\alpha<2$, \citet[Theorem~9]{BBS2007} showed
that, for each fixed $k$, the number $N_k(n)$ of blocks of size $k$ of the allelic partition
satisfies $n^{\alpha-2}N_k(n)\to\tfrac\theta2\alpha(\alpha-1)^2\Gamma(k+\alpha-2)/k!$ in
probability; their mutation rate, per unit length of lineage, is our $\theta/2$. In
particular $N_2(n)/N_1(n)\to(\alpha-1)/2$ in probability. With $\sigma=2-\alpha$, the
limiting proportions $(2-\alpha)\Gamma(k+\alpha-2)/(\Gamma(\alpha-1)k!)$ that their
Remark~13 associates with this theorem equal $\sigma\rf{1-\sigma}{k-1}/k!$, which involves
the Gibbs weights of type $2-\alpha$. The edge index converges to a smaller value:
\[
\sigma_\star(\alpha)-(2-\alpha)=-\frac{\alpha(2-\alpha)^2}{1+3\alpha-\alpha^2}<0
\qquad(1<\alpha<2),
\]
since $1+3\alpha-\alpha^2=3+(\alpha-1)(2-\alpha)\ge3$ there. Thus the small blocks of a
typical large sample involve the Gibbs weights of type $2-\alpha$, while the conditional law
of $\Pi_n$ given $K_n=n-2$, an event of vanishing probability by
Lemma~\ref{lem:rare}, is governed by an edge index whose limit $\sigma_\star(\alpha)$ is
smaller than $2-\alpha$.

Table~\ref{tab:index} lists exact values of $\sigma^{(n)}$, rounded, together with the
limits given by Proposition~\ref{prop:limits}.

\begin{table}
\caption{Edge index $\sigma^{(n)}$, rounded to four decimals, for four measures $\Lambda$
and two values of $\theta$, with its limit as $n\to\infty$ from
Proposition~\ref{prop:limits}. $\Beta(1/2,3/2)$ and $\Beta(3/2,1/2)$ are the cases
$\alpha=3/2$ and $\alpha=1/2$ of $\Beta(2-\alpha,\alpha)$, and BS denotes the
Bolthausen--Sznitman coalescent.}\label{tab:index}
\begin{tabular}{@{}lccccl@{}}
\toprule
$\Lambda$, $\theta$ & $n=4$ & $n=10$ & $n=100$ & $n=1000$ & limit\\
\midrule
$\delta_0+\delta_1$, $\theta=1$ & 0.3846 & 0.0553 & 0.0004 & 0.0000 & $0$\\
$\Beta(1/2,3/2)$, $\theta=1$ & 0.4118 & 0.4059 & 0.3946 & 0.3883 & $5/13\approx0.3846$\\
$\Beta(1/2,3/2)$, $\theta=5$ & 0.6000 & 0.5318 & 0.4400 & 0.4037 & $5/13$\\
BS, $\theta=1$ & 0.6923 & 0.6917 & 0.7070 & 0.7129 & $5/7\approx0.7143$\\
BS, $\theta=5$ & 0.8400 & 0.8218 & 0.8166 & 0.8177 & $9/11\approx0.8182$\\
$\Beta(3/2,1/2)$, $\theta=1$ & 0.8824 & 0.8908 & 0.9391 & 0.9715 & $1$\\
\bottomrule
\end{tabular}
\end{table}

\section{Open problems}\label{sec:open}

Three questions remain open. First, whether some $\Lambda$ and some $\theta>0$ satisfy the
edge identity \eqref{eq:edge-identity} for one $\sigma\in(0,1)$ at every $n\ge4$, that is,
give an edge index that does not depend on $n$; by Proposition~\ref{prop:limits}(a) and (d),
whatever $\theta$, such a measure has no atom at $0$
and no density $f$ with $f(x)\sim Cx^{1-\alpha}$ as $x\downarrow0$ for some $C>0$ and
$0<\alpha<1$, and by Proposition~\ref{prop:sharp}(b) and the scaling of
Section~\ref{sec:model} it is not a multiple of the Lebesgue measure. Second, whether
\[
\PP(n\in B\mid\Pi_{n-1}=\pi)\ge|B|\,\PP(n\in\{j\}\mid\Pi_{n-1}=\pi)
\]
for every $\Lambda$, every $\theta>0$, every $n\ge3$, every partition $\pi$ of $[n-1]$ with
$\PP(\Pi_{n-1}=\pi)>0$, every block $B$ of $\pi$ and every singleton block $\{j\}\ne B$ of
$\pi$.
Under the Ewens formula equality holds, and Theorem~\ref{thm:edge}(b) is the case where $B$
has two elements and all other blocks are singletons. Between two blocks of
size at least two the analogous inequality fails: for $\Lambda=\delta_{1/2}$, $\theta=1$ and
$\pi=\{\{1,2,3\},\{4,5\}\}$, given $\Pi_5=\pi$ the probability that $6$ joins $\{1,2,3\}$ is
$3633/2554<3/2$ times the probability that it joins $\{4,5\}$. Third, whether Lemma~\ref{lem:certificate},
or the exclusion of types $\sigma\in(0,1)$ at a fixed $\theta$, has a proof that does not
go through the coefficients of $Q$.

\appendix
\section{The certificate polynomials}\label{app:certificate}

Throughout, $h_0=1$ and $a$, $b$, $c$ are the coordinates of Lemma~\ref{lem:certificate},
so that $\sigma=a/(1+a)$, $\theta=b$ and $h_1=c/(1+c)$; we write
\[
x=b(1+c),\qquad y=1+3c .
\]
The formulas below come from the same exact computation as \eqref{eq:certificate}. The
proof of Theorem~\ref{thm:rigidity} uses only that the slopes are negative at the moments
of a measure, which is Lemma~\ref{lem:affine}.

\paragraph{Slopes and rates}
The slopes of Lemma~\ref{lem:affine}, evaluated along $h_2=\eta_2$ and $h_3=\eta_3$, are
\begin{gather*}
\kappa_4=-\frac{\ell_4}{(1+a)(1+c)},\qquad \kappa_5=-\frac{4\ell_5}{(1+a)(1+c)\ell_4},\\
\kappa_6=-\frac{6\ell_6}{(1+a)(1+c)\ell_4\ell_5},
\end{gather*}
with
\[
\ell_4=(3+2a)y+3x,\qquad \ell_5=2(a+2)(a+3)y^2+3(4a+7)xy+9x^2,
\]
and $\ell_6$ a polynomial of total degree $8$ in $(a,b,c)$ with $75$ terms, all with positive
coefficients. Formula \eqref{eq:H2} becomes $\eta_2=y(x+y)/((1+c)\ell_4)$. The total rates at
the moments $(1,h_1,\eta_2,\eta_3,\eta_4)$ are
\begin{gather*}
\rho_2=1+b,\qquad\rho_3=\frac{2y+3x}{2(1+c)},\qquad\rho_4=\frac{2F_4}{(1+c)\ell_4},\\
\rho_5=\frac{F_8}{2(1+c)\ell_4\ell_5},\qquad\rho_6=\frac{F_{16}}{3(1+c)\ell_4\ell_5\ell_6},
\end{gather*}
where $F_4=(a+3)y^2+2(a+3)xy+3x^2$, and $F_8$, $F_{16}$ are the polynomials of
Lemma~\ref{lem:certificate}. Thus all slopes are negative and all rates positive on
$(0,\infty)^3$. Since $1+b+3c+bc=x+y$ and $2+3b+6c+3bc=2y+3x$, the denominator of
\eqref{eq:certificate} is $2+a=(1+a)(2-\sigma)$ times the product of the polynomials $1+b$,
$2y+3x$, $F_4$, $F_8$ and $F_{16}$ that appear in the numerators of $\rho_2,\dots,\rho_6$.
The polynomial $Q$ has degree $8$ in $a$ and $9$ in each of $b$ and $c$.

\paragraph{The data file}
The supplementary file \texttt{certificate-n6.json} lists $Q$, $F_4$, $F_8$, $F_{16}$,
$\ell_4$, $\ell_5$ and $\ell_6$ (under the keys \texttt{k4}, \texttt{k5}, \texttt{k6}), each as a list of entries $[[i,j,k],\,q]$ meaning
$q\,a^ib^jc^k$, together with identity \eqref{eq:certificate}, the formulas for the rates
above, and the assertion that every coefficient is a positive integer.

\section*{Funding}
This work was supported in part by the Conselho Nacional de Desenvolvimento
Cient\'ifico e Tecnol\'ogico (CNPq) [grant numbers 301135/2026-6 (R.O.), 310800/2026-9
and 307626/2022-9 (A.M.S.M.)]; and the Coordena\c{c}\~ao de Aperfei\c{c}oamento de
Pessoal de N\'ivel Superior (CAPES) [Finance Code 001 (R.R.A.)]. The funding sources
had no involvement in the study design, in the analysis and interpretation of the
results, in the writing of the article, or in the decision to submit it for publication.

\section*{Declaration of competing interest}
The authors declare that they have no known competing financial interests or personal
relationships that could have appeared to influence the work reported in this paper.

\section*{Data availability}
The polynomials of Lemma~\ref{lem:certificate} and Appendix~\ref{app:certificate} are
supplied as the supplementary file \texttt{certificate-n6.json}, as lists of exponent
vectors with integer coefficients: $Q$ ($530$ monomials, total degree $18$), $F_4$, $F_8$
and $F_{16}$ ($15$, $70$ and $396$ monomials), and $\ell_4$, $\ell_5$ and $\ell_6$ ($6$, $18$ and
$75$ monomials). The programs that derive these polynomials and check the identities and
exact values of Sections~\ref{sec:edge}--\ref{sec:limits} are in the folder
\texttt{arxiv-v2} of the public repository \repourl~\citep{GibbsRepo2026}; the programs and
certificates of the first preprint version remain in the same repository.

\section*{Declaration of generative AI and AI-assisted technologies in the manuscript preparation process}
During the preparation of this work the authors used Claude (Anthropic) to develop and
test code, structure the manuscript and revise the text. After using this tool, the authors
reviewed and edited the content as needed and take full responsibility for the content of
the published article.

\section*{Acknowledgements}
The authors thank the members of the Laborat\'orio de F\'isica Te\'orica e
Computacional and of the Departamento de Estat\'istica of the Universidade Federal de
Pernambuco for discussions.

\end{document}